\documentclass[11pt,a4paper]{article}
\usepackage{amsmath}
\usepackage{amsthm}
\usepackage{amsfonts}
\usepackage{latexsym}
\usepackage{geometry}
\newcounter{alphthm}
\newtheorem{thm}{Theorem}[section]

\newtheorem{lem}[thm]{Lemma}

\theoremstyle{definition}

\newcommand{\be}{\begin{equation}}
\newcommand{\ee}{\end{equation}}
\newcommand{\ben}{\begin{enumerate}}
\newcommand{\een}{\end{enumerate}}

\def\beq{\begin{equation}}
\def\eeq{\end{equation}}

\title{On Projectively Flat Spherically Symmetric Finsler Metrics}
\author{B. Najafi}

\begin{document}
 \maketitle

\begin{abstract}
The class of spherically symmetric Finsler metrics is studied  and  locally dually flat and locally projectivelly flat spherically symmetric Finsler metrics is classified.
\end{abstract}

\section{Introduction}
Let $F$ be a Finsler metric defined on a convex  domain $\Omega\subset \mathbb{R}^n$
and  be invariant under any rotation in $\mathbb{R}^n$. Then, $F$ is called spherically symmetric. In \cite{Zh}, by solving the equation of Killing fields, Zhou showed that there exists a positive function $\phi$  so that $F$ can be written as $F = \phi(r,u,v)$  where
\[
r = |x|, \ \  \ u = |y|, \ \ \ v = <x,y>, \ \  \ s =\frac{<x,y>}{|y|}
\]
and  $|.|$ and $<,>$ denote the Euclidean norm and inner product in
$\mathbb{R}^n$, respectively.  Many well-known classical Finsler metrics such as Funk metric and Berwald metric are spherically symmetric. Having a nice symmetry makes the class of spherically symmetric Finsler metrics very important both in mathematics and applications \cite{MST}\cite{MZ}.

A Finsler metric is said to be  locally projectively flat if at any point there is a local coordinate system  in which the geodesics are straight lines as point sets. Recently, motivated by Hilbert's Fourth problem relating to classify the projectively flat  Finsler metrics in $\mathbb{R}^n$, Zhou completely classified projectively flat  spherically  symmetric Finsler metrics \cite{Zh}\cite{GLM}. According to Rapcs\'{a}k  Lemma, a Finsler metric $F$ on an open subset $\mathcal{U}\in R^n$ is  projectively flat on
$\mathcal{U}$ if and only if it satisfies $F_{x^ky^l}y^k=F_{x^l}$.

In \cite{amna}, Amari-Nagaoka introduced the notion of dually flat Riemannian metrics when they study the information geometry on Riemannian manifolds.  In Finsler geometry, Shen extends the notion of locally dually flatness for Finsler metrics \cite{shenLec}.  A Finsler metric $F = F(x, y)$ on an open subset $U\subset \mathbb{R}^n$ is dually flat if and only if it satisfies
\[
[F^2]_{x^my^l}y^m=2[F^2]_{x^m}.
\]

In this paper, we characterize locally dually flat spherically  symmetric Finsler metrics and give complete classification of projectively flat metrics among them. More precisely, we have the following.

\begin{thm} \label{main theorem}
Let $F = \phi(r,u,v)$ be a spherically symmetric Finsler metric on a convex  domain $\Omega\subset \mathbb{R}^n$. Then $F$ is projectively flat and locally dually flat if and only if
\be
\phi(r,u,v)=\frac{\sqrt{(k-c^2r^2)u^2+c^2v^2}+cv}{k-c^2r^2},\label{TSN}
\ee
where $k$ and $c$ are constants. More precisely, every  projectively flat and locally dually flat spherically symmetric Finsler metrics is a deformation of Funk metric.
\end{thm}

\bigskip
It is remarkable that, the Funk metric $F$ is  defined on the standard unite ball $B^n$ as follows
\[
F(x,y):=\frac{\sqrt{|y|^2-(|x|^2|y|^2-\langle x,y\rangle^2)}+\langle
 x,y\rangle}{1-|x|^2}.
\]
Thus  the Funk metric is a special case of the metric defined by (\ref{TSN}) with  $k=c=1$.

\section{Preliminary}

Given a Finsler manifold $(M,F)$, then a global vector field ${\bf G}$ is induced by $F$ on slit tangent bundle  $TM_0=TM-\{0\}$, which in a standard coordinate $(x^i,y^i)$ for $TM_0$ is given by
\[
{\bf G}=y^i {{\partial} \over {\partial x^i}}-2G^i(x,y){{\partial} \over {\partial y^i}},
\]
where $G^i=G^i(x, y)$ are called spray coefficients and given by
\[
G^i=\frac{1}{4}g^{il}\Big[\frac{\partial^2F^2}{\partial x^k \partial y^l}y^k-\frac{\partial F^2}{\partial x^l}\Big].
\]
${\bf G}$ is called the  spray associated  to $F$.

\bigskip

A Finsler metric $F=F(x, y)$ is called  locally projectively flat if at any point there is a local coordinate system  in which the geodesics are straight lines as point sets. It is known that a Finsler metric  $F(x,y)$ on an open domain $ U\subset \mathbb{R}^n$   is  locally projectively flat  if and only if its geodesic coefficients $G^i$ are in the form
\[
G^i= Py^i,
\]
where $P: TU = U\times \mathbb{R}^n \to \mathbb{R}$ is positively homogeneous with degree one,
$P(x,  \lambda y) = \lambda P(x, y)$, $\lambda >0$. We call $P(x, y)$ the  projective factor of $F(x, y)$.

\bigskip

A Finsler metric $F=F(x,y)$ on a manifold $M$ is said to be locally dually flat if at any point there is a  coordinate system $(x^i)$ in which the spray coefficients are in the following form
\[
G^i = -\frac{1}{2}g^{ij}H_{y^j},
\]
where $H=H(x, y)$ is a $C^\infty$ homogeneous scalar function on $TM_{0}$. Such a coordinate system is called an adapted coordinate system. In \cite{shenLec}, Shen proved that  the Finsler metric $F$ on an open subset $U\subset \mathbb{R}^n$ is dually flat if and
only if it satisfies $(F^2)_{x^ky^l}y^k=2(F^2)_{x^l}$. In this case, $H =-\frac{1}{6}[F^2]_{x^m}y^m$.

\bigskip

A Finsler metric $F$ on  a domain $\Omega\subseteq \mathbb{R}^n$ is called spherically symmetric if it is invariant under any rotation in $\mathbb{R}^n$. According to the equation of Killing fields,  there exists a positive function $\phi$ depending on two variables so that $F$ can be written as
\[
F=|y|\phi\Bigg(|x|, \frac{\langle x,y\rangle}{|y|}\Bigg),
\]
where $x$ is a point in the domain $\Omega$, $y$ is a tangent vector at the point $x$ and $\langle, \rangle$, \  $|\cdot|$ are standard inner product and norm in Euclidean space, respectively.

\medskip

\begin{lem} {\rm (\cite{Zh})}
\emph{A Finsler metric $F$ on  a convex domain $\Omega\subseteq \mathbb{R}^n$ is spherically symmetric if and only if
there exists a positive function $\phi=\phi(r,u,v)$, such that $F(x,y)=\phi(|x|,|y|,\langle x,y\rangle)$, where $|x|=\sqrt{\Sigma^n_{i=1}(x^i)^2}$,
$|y|=\sqrt{\Sigma^n_{i=1}(y^i)^2}$ and $\langle
x,y\rangle=\Sigma^n_{i=1}x^iy^i$.}
\end{lem}

%--------------------------------------------------------------------------------------------------------------------
\section{Proof of Theorem \ref{main theorem}}

A Finsler metric $F = F(x, y)$ on an open subset $U\subset \mathbb{R}^n$ is dually flat if and only if it satisfies
\be
[F^2]_{x^my^l}y^m=2[F^2]_{x^m}. \label{2}
\ee
In \cite{CSZ}, X. Cheng, Z. Shen and Y. Zhou studied  and characterized projectively and locally dually flat Finsler metrics on a convex  domain $\Omega\subset \mathbb{R}^n$ and found the following PDEs.
\medskip
\begin{thm}\label{theorem1}
Let $F = F(x, y)$ be a Finsler metric on an open subset $U\subset \mathbb{R}^n$. Then
$F$ is dually flat and projectively flat on $U$ if and only if it satisfies
\be
F_{x^k}=c F F_{y^k}\label{0}
\ee
where $c$ is a constant.
\end{thm}

\bigskip

For a spherically symmetric  Finsler metric in
$\mathbb{R}^n$, we have the following.
\medskip
\begin{thm}\label{thm2}
Let $F = \phi(r, u, v)$ be a  spherically symmetric  Finsler metric in
$\mathbb{R}^n$. Then $F$ is locally dually flat if and only if the following holds
\begin{eqnarray}
s(\psi_r\psi_s+\psi\psi_{rs})+r(\psi_s^2+\psi\psi_{ss})-2\psi\psi_r=0,\label{3}
\end{eqnarray}
where
\[
s:=\frac{v}{u}\ \ \  \ \textrm{and}  \ \ \  \ \psi(r,s):=\phi(r,1,\frac{v}{u}).
\]
\end{thm}
\noindent
\textbf{\textsl{Proof}.}
By direct computations, we have
\begin{eqnarray}
&&F_{x^l}=\frac{u}{r}\psi_r x^l+\psi_s y^l,\label{5}\\
&&F_{x^k}y^k=u^2\Big[\frac{ s}{r}\psi_r+\psi_s\Big],\label{6}\\
&&F_{y^l}=\psi_s x^l+\frac{1}{u}\Big[\psi-s\psi_s\Big]y^l,\label{7}\\
&&F_{x^ky^l}y^k=u\Big[\frac{s}{r}\psi_{rs}+\psi_{ss}\Big]x^l+\Big[\frac{s}{r}\psi_r-\frac{s^2}{r}\psi_{rs}-\psi_{ss}s
+\phi_{s}\Big]y^l.\label{8}
\end{eqnarray}
By (\ref{2}),  $F$ is locally dually flat if and only if
\begin{eqnarray}
F_{y^l}F_{x^k}y^k+FF_{x^ky^l}y^k-2FF_{x^l}=0.\label{9}
\end{eqnarray}
Plugging (\ref{5}), (\ref{6}), (\ref{7}) and (\ref{8}) in (\ref{9}) imply (\ref{2}). This completes the proof.
\qed

\bigskip

In \cite{Zh},  L. Zhou gave the following classification of projectively spherically symmetric Finsler metric.
\medskip

\begin{thm}\label{theorem2}
Suppose that $F$ is a spherically symmetric Finsler metric on a convex
domain  $\Omega\subset \mathbb{R}^n$, $F$ is projectively flat if and only if there exist smooth functions $f=f(t) > 0$
and $g=g(r)$ such that
\begin{eqnarray}
\phi(r,u,v)=\int f(\frac{v^2}{u^2}-r^2)du+g(r)v, \label{Eq1}
\end{eqnarray}
where $F(x,y)=\phi(|x|,|y|,\langle x,y\rangle)$.
\end{thm}

\bigskip

Now, we are going to prove the main result.

\bigskip
\noindent
{\bf Proof of Theorem \ref{main theorem}:} Using the identities
\[
F_{x^l}=\frac{\phi_r}{r} x^l+\phi_v y^l, \ \ \  and  \ \ \ F_{y^l}=\phi_v x^l+\frac{\phi_u}{u}y^l
\]
and by Theorem \ref{theorem1}, we get
\begin{eqnarray}
&&\frac{\phi_r}{r}=c\phi\phi_v,\label{Eq11}\\
&&\phi_v=\frac{c\phi\phi_u}{u}.\label{Eq12}
\end{eqnarray}
 By  (\ref{Eq1}) we have
\be
\phi_u(r,u,v)=f(\frac{v^2}{u^2}-r^2).\label{e1}
\ee
Plugging (\ref{e1}) into  (\ref{Eq12})  yields
  \begin{eqnarray}
 \phi_v(r,u,v)=\frac{c}{u}\phi(r,u,v)f\Big(\frac{v^2}{u^2}-r^2\Big)\label{Eq3}.
 \end{eqnarray}
Using 1-homogeneity of $\phi$ with respect to $(u,v)$ and Euler's theorem, we have $\phi=\phi_u u+\phi_v v$. Thus, from (\ref{Eq3}) we conclude that
 \begin{eqnarray}
 \phi(r,u,v)=\frac{f(\frac{v^2}{u^2}-r^2)u^2}{u-c f(\frac{v^2}{u^2}-r^2)v}\label{Eq4}.
 \end{eqnarray}
 Therefore, it suffices to find explicit formula of $f$.  Taking derivative of  (\ref{Eq4}) with respect to $r$ implies that
 \begin{eqnarray}
 \frac{\phi_r(r,u,v)}{r}=\frac{-2f'(\frac{v^2}{u^2}-r^2)u^3}{\Big[u-c f(\frac{v^2}{u^2}-r^2)v\Big]^2}\label{Eq5}
 \end{eqnarray}
Substituting  (\ref{Eq11}) and (\ref{Eq3}) into (\ref{Eq5}), we obtain the following ODE on $f$:
 \be
 2f'+c^2f^3=0\label{Eq6}
 \ee
 Solving (\ref{Eq6}), we have
\[
f(t)=\frac{1}{\sqrt{c^2t+k}},
\]
where $k$ is a constant. The proof follows from (\ref{Eq4}).
\qed

\bigskip

\noindent
Behzad Najafi\\
Department of Mathematics and Computer Sciences\\
Amirkabir University\\
Tehran. Iran\\
Email:\ behzad.najafi@aut.ac.ir

\end{document}